\documentclass[12pt]{amsart}

\usepackage{amssymb}

\usepackage{epsfig,color}

\newtheorem{thm}{Theorem}
\newtheorem{prop}{Proposition}
\newtheorem{defin}{Definition}

\newtheorem{corol}{Corollary}
\newtheorem{lem}{Lemma}
\newtheorem{conj}{Conjecture}
\newtheorem{rem}{Remark}

\newtheorem{exa}{Example}
\newtheorem{quest}{Question}

\renewcommand{\c}{{\mathcal C}}

\newcommand{\II}{{\mathbb I}}

\newcommand{\JJ}{{\mathbb J}}

\newcommand{\n}{{\mathfrak n}}\newcommand{\N}{{\mathbb N}}

\newcommand{\QQ}{{\mathbb Q}}

\newcommand{\R}{{\mathbb R}}\newcommand{\RR}{{\mathbb R}^2}

\newcommand{\tf}{{\mathbb T}}

\newcommand{\T}{{\mathbb T}}

\newcommand{\Z}{{\mathbb Z}}

\newcommand{\bo}{\partial} 

\newcommand{\const}{\mbox{const}} 

\newcommand{\hg}{\hat {g}}

\newcommand{\tA}{\tilde{A}}

\newcommand{\tS}{\tilde{S}}

\newcommand{\al}{\alpha}
\newcommand{\be}{\beta}
\newcommand{\ga}{\gamma}\newcommand{\Ga}{\Gamma}
\newcommand{\de}{\delta}

\newcommand{\ka}{\kappa}

\newcommand{\Om}{\Omega}
\newcommand{\si}{\sigma}
\newcommand{\vp}{\varphi}
\newcommand{\het}{\theta}

\begin{document}

\bibliographystyle{plain}

\title[Floating and billiard]
{Capillary floating and the billiard ball problem}

\author{Eugene Gutkin}

\address{Department of Mathematics, Nicolaus Copernicus University, Chopina 12/18, Torun 87-100, Poland;
Institute of Mathematics of Polish Academy of Sciences,
Sniadeckich 8, Warsaw 00-956, Poland}
\email{gutkin@mat.umk.pl,gutkin@impan.pl}

\keywords{floating in neutral equilibrium, surface tension,
contact angle, billiard ball problem, invariant curves, constant
angle caustics, domains of constant width}

\subjclass{76B45,76D45,52A10,37A45,37E10,37E40}

\date{\today}

\begin{abstract}
We establish a connection between capillary floating in neutral
equilibrium and the billiard ball problem. This allows us to
reduce the question of floating in neutral equilibrium at any
orientation with a prescribed contact angle for infinite
homogeneous cylinders to a question about billiard caustics for
their orthogonal cross-sections. We solve the billiard problem. As
an application, we characterize the possible contact angles and
exhibit an infinite family of real analytic non-round cylinders
that float in neutral equilibrium at any orientation with constant
contact angles.
\end{abstract}

\maketitle

\tableofcontents

\section{Introduction: Floating in neutral equilibrium and the billiard ball problem}       \label{intro}
The mathematical theory of capillarity goes back to 1806. In his
famous treatise on celestial mechanics~~\cite{La1806} Laplace
discussed a broad range of problems related to surface tension at
fluid interfaces, among them a theory of capillary floating. One
of the major open problems in this subject is to determine
configurations at which a particular body will float on a liquid
surface. In~~\cite{La1806} Laplace characterized some special
cases of capillary floating which was an astonishing achievement
for his time. There are several physical phenomena that need to be
taken into account: The mass distribution in the body, the
gravity, the surface tension, etc. This leads to a highly
nonlinear free boundary problem with nonlinear boundary
conditions. Further specializing and simplifying the physical
assumptions, we arrive at the concept of {\em floating in neutral
equilibrium with zero gravity}. See \cite{F09i}. In what follows
we will simply speak of floating in neutral equilibrium. The
relevant mathematical conditions involve geometry of the body
surface, its space orientation, and the {\em contact angle}
between the body and the liquid surface.

\begin{figure}[htbp]
\input{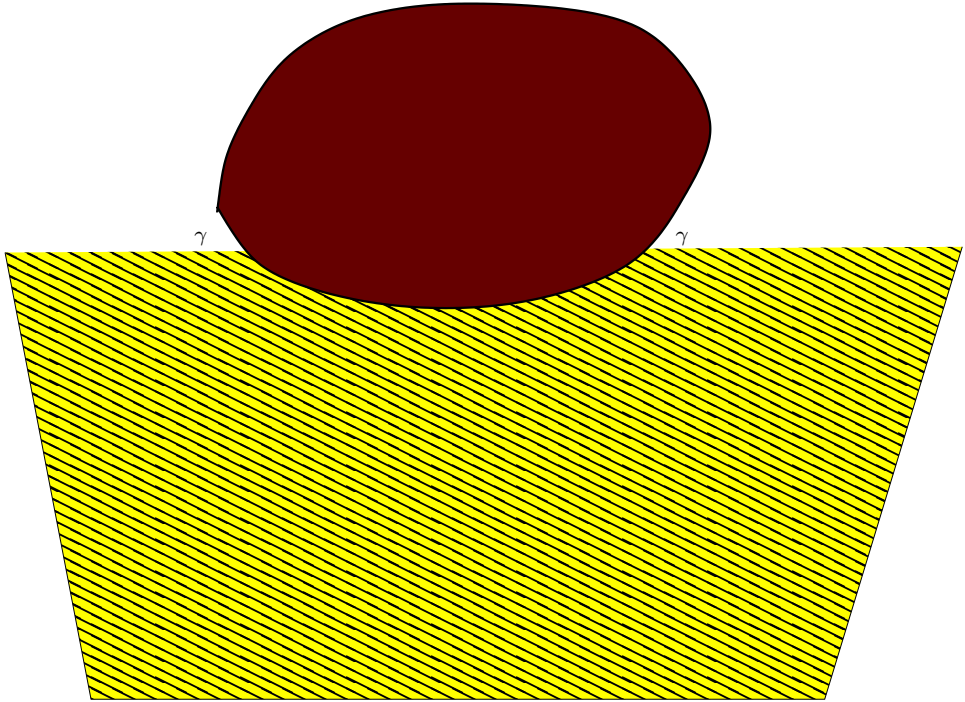} 
\begin{center}
\caption{Floating in neutral equilibrium; $\ga$ is the contact
angle.} \label{neut_equi_fig}
\end{center}
\end{figure}

On the outset, the problem of body floating is three dimensional.
In the special case when the body is an infinite homogeneous
cylinder, it reduces to a two dimensional problem involving the
cross-section of the cylinder. We will assume that it is a
bounded, convex, planar domain, say $\Om\subset\RR$. The
three-dimensional floating problem for the cylinder translates
into properties of $\Om$ which are naturally interpreted as the
conditions for capillary floating in two dimensions. The essential
requirements for neutral equilibrium are then imposed by the basic
laws of physics. In what follows we will be mostly concerned with
the two-dimensional floating in neutral equilibrium. This problem
can be reformulated in terms of the geometry of $\Om$.
Figure~~\ref{neut_equi_fig} illustrates the concept of
two-dimensional floating in neutral equilibrium.

The first work connecting two-dimensional capillary floating with
convex geometry appears to be~~\cite{RMBC}. The main result
in~~\cite{RMBC} says that any regular, convex, bounded, planar
domain will float at a given contact angle in at least four
distinct orientations. The proof crucially uses the four vertex
theorem. The problem was further studied by R. Finn.
In~~\cite{F09i} he showed that the mathematical assumptions
in~~\cite{RMBC} follow from basic physical laws. Finn then asked
which convex, smooth domains $\Om$ can float at a prescribed
contact angle $\ga$ in every rotational orientation. He pointed
out that for $\ga=\pi/2$ this happens if and only if $\Om$ is a
domain of constant width. See \cite{F09i,F09ii} for this and
related material.

The geometry of a convex planar domain $\Om$ is crucial for the
billiard ball problem championed by G.D. Birkhoff in the early
20th century~~\cite{Bi17}. It can also be viewed as a highly
specialized and simplified case of a physical situation. The
billiard ball is a point that travels with the unit speed inside
$\Om$ reflecting at the boundary $\bo\Om$ according to the law of
equal angles. Disregarding the motion of the ball between
collisions with the boundary, we reduce the billiard ball problem
to the study of the {\em billiard map} on $\Om$. Invariant curves
of this map provide a crucial insight into the billiard dynamics.
Let $s$ be the arc length variable on $\bo\Om$, and let $\het$ be
the outgoing angle. See Figure~~\ref{bill_map_fig}. Beginning with
Birkhoff~~\cite{Bi17}, invariant curves of the form $\het=h(s)$
played an important role in the literature on billiard dynamics.

\medskip

The functions $\het=h(s)$ that yield invariant curves have been
extensively studied~~\cite{Laz73,Dou,In88,GK}. The present work is
based on the following observation relating the floating problem
and the billiard ball problem for the planar domain $\Om$. The
cylinder with the cross-section $\Om$ floats in neutral
equilibrium at any orientation with the contact angle $\ga$ if and
only if the billiard table $\Om$ admits the invariant curve
$\het=h(s)$ with the constant function $h(s)\equiv\pi-\ga$. For
the reasons that we will explain in section~~\ref{birkhoff}, we
call these invariant curves the {\em constant angle caustics} for
$\Om$. The floating problem thus becomes the following billiard
problem: i) Find the regular, convex billiard tables that admit
constant angle caustics; ii) Determine the corresponding angles.
This work provides a fair amount of information on this subject.
Before describing our results, we will further elaborate on the
capillary floating in three dimensions.

\medskip

The conjecture that the round ball is the only body to float  in
equilibrium at any orientation is usually ascribed to S. Ulam. See
\cite{Au38,Scott}. Various authors have mathematically
reformulated this question in different, albeit related ways. The
interpretation of Finn et al takes the tension of the liquid
surface   into account~~\cite{F09i,FS09,FV09}, while the other
interpretations disregard it~~\cite{Va}. Whatever the
interpretation, the conditions of floating in neutral equilibrium
are much more restrictive for three dimensions than in the special
case of two dimensions. Thus, Finn and Sloss~~\cite{FS09} show
that the only three-dimensional body to float in neutral
equilibrium in any orientation at a constant contact angle is the
round ball. Using a different interpretation of the concept of
floating, P. Varkony~~\cite{Va} finds a counterexample to the Ulam
conjecture. It is clear that the concept of floating in neutral
equilibrium generates challenging questions about the geometry of
surfaces in $\R^3$. In the present work we relate the
two-dimensional floating to the geometry of convex bounded domains
in $\RR$.

\medskip

\begin{figure}[htbp]
\input{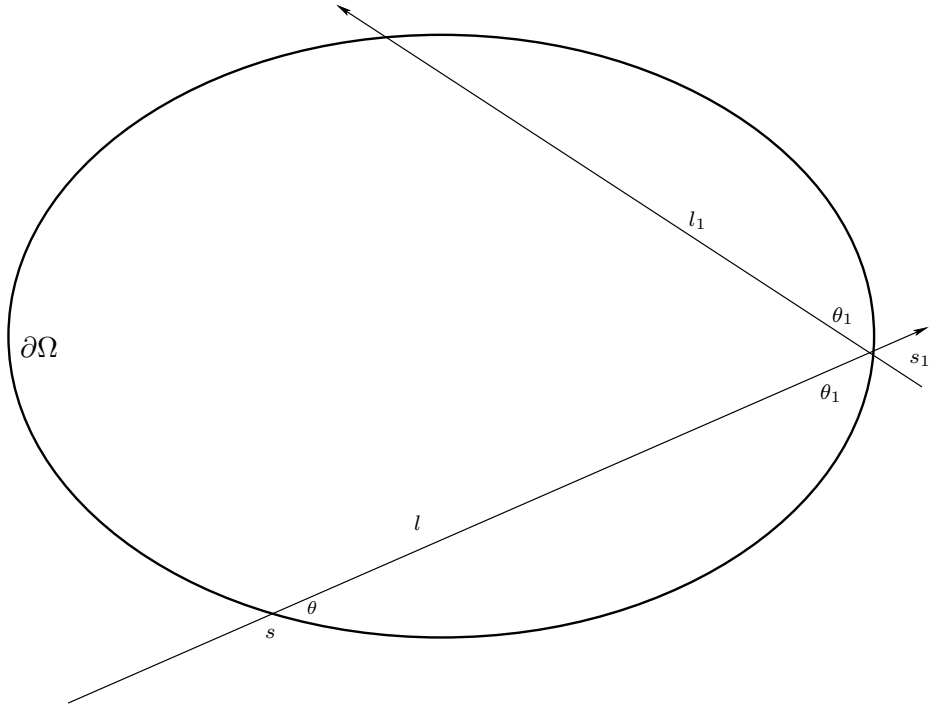} 
\begin{center}
\caption{The billiard phase space as a space of rays intersecting
the billiard table.} \label{light_rays_fig}
\end{center}
\end{figure}

We will now briefly describe our results and the structure of the
paper. In section~~\ref{birkhoff} we review the concept of the
{\em billiard map}. Let $\Om\subset\RR$ be a bounded, strictly
convex domain with the smooth boundary $\bo\Om$. The {\em phase
space} $Z=Z(\Om)$ for the billiard map on $\Om$ consists of
rays\footnote{I. e., oriented straight lines.} intersecting $\Om$.
Let $0\le\het,\het_1\le\pi$ be the two angles that a ray $l\in Z$
forms at the points of intersection $s,s_1\in\bo\Om$. See
Figure~~\ref{light_rays_fig}. Let $F:Z\to Z$ be the billiard map.
The ray $l_1=F(l)$ is obtained by reflecting $l$ at $s_1$ about
$\bo\Om$, as if $\bo\Om$ was a perfect mirror. The domain $\Om$
floats in neutral equilibrium at any orientation with a constant
contact angle if and only if there exists $0<\de<\pi$ such that
for any $l\in Z$ satisfying $\het(l)=\de$ we have $\het_1(l)=\de$.

In section~~\ref{birkhoff} we study this geometric condition from
the viewpoint of the billiard map. Let $\rho(s)$ be the radius of
curvature for $\bo\Om$; let $c_k,k\in\Z,$ be its Fourier
coefficients. Note that $\bo\Om$ is circular if and only if
$c_k=0$ for all nonzero $k$. Theorem~~\ref{const_caust_thm} says
that a noncircular domain $\Om$ has the above property if and only
if the following two conditions hold: i) There exists $n>1$ such
that the pair $n,\de$ satisfies the trigonometric
equation~~\eqref{sin_eq}; ii) The coefficients $c_k$ vanish if the
pair $k,\de$ does not satisfy equation~~\eqref{sin_eq}.

In section~~\ref{reduce} and section~~\ref{impli} we study
equations~~\eqref{sin_eq} and obtain several applications. The
value $\de=\pi/2$ is special in that the pair $n,\pi/2$ satisfies
equation~~\eqref{sin_eq} with any odd $n$. The corresponding
regions $\Om$ are the domains of constant width; we briefly review
their geometry in section~~\ref{const_width}. The symmetry
$\de'=\pi-\de$ allows us to reduce the study of solutions of
equations~~\eqref{sin_eq} to the range $0<\de<\pi/2$. Restricted
to this interval, equations~~\eqref{sin_eq} are equivalent to
$\tan n x=n\tan x$. In section~~\ref{trig_qual} we obtain fairly
detailed qualitative information about solutions of these
equations. Let $B_n\subset(0,\pi/2)$ denote the set of solutions.
We show that $B_n$ has roughly $n/2$ elements; it is
$(\pi/n)$-dense in the interval $(0,\pi/2)$. For every $n>3$ we
exhibit  a one-parameter family $\Om_{n,\tau}$ of noncircular,
real analytic domains that float in every orientation at the
contact angles $\ga\in B_n$ and $\ga\in \pi-B_n$. See
equations~~\eqref{table_eq1},~~\eqref{table_eq2} and
Corollary~~\ref{float_cor}.

The classification of domains that float in neutral equilibrium at
constant contact angles hinges on the information about the
solutions to $\tan n x=n\tan x$. In particular, we need to know
whether the sets $B_n\cap B_m$ can have nonempty intersections for
$m\ne n$. In section~~\ref{trig_quant} and
section~~\ref{numb_theor} we reduce these questions to a study of
roots of an infinite chain of polynomials that are closely related
to Chebyshev polynomials. This reveals a number-theoretic aspect
of capillary floating. Let $S_n,n\ge 1$, be the polynomials. In
section~~\ref{frac_lin_sub} we study the roots  of $S_n$ and
obtain some information about them. However, the question whether
$S_m,S_n$ have nontrivial common roots for $m\ne n$ remains
unresolved. It is the lack of this information that prevented the
author from publishing his findings immediately after the 1993
PennState Dynamics Workshop~~\cite{Gut93}. The book~~\cite{Ta95}
contains a brief report on these findings.

Based on substantial partial evidence, we formulate three
conjectures about the roots of $S_n$.
Conjectures~~\ref{disjoint_roots_conj} and~~\ref{dist_tan_conj}
are equivalent.
In section~~\ref{cond_impli}, assuming that these conjectures
hold, we derive consequences for the billiard and for the floating
problem.
Theorems~~\ref{cond_const_caust_thm},~~\ref{irrat_rot_thm},~~\ref{conseq_thm}
completely describe billiard tables with constant angle caustics.
Theorem~~\ref{float_thm} gives a classification of regular planar
domains that float in neutral equilibrium in any orientation at
constant contact angles.

\medskip

The present work can be viewed as one of many examples of fruitful
relationships between the billiard and other mathematical
subjects. We refer the reader to
\cite{Ka05,KH,GR09,Ta95,Gut97,Ta05} for other examples of this
nature. The billiard framework offers a variety of open problems
that often bear on fundamental and elementary mathematical
concepts \cite{Gut03}. The author hopes that the present work will
help to advertise  the subject in the mathematical fluid mechanics
community. The author is grateful to Bob Finn for bringing the
subject of capillary floating to his attention and for making
several comments on the present work. It was partially supported
by the MNiSzW grant N N201 384834.

\section{The Birkhoff billiard: General caustics versus constant angle caustics}    \label{birkhoff}
The billiard in the sense of G.D. Birkhoff plays on a compact,
convex domain $\Om\subset\RR$. We will assume that the boundary
$\bo\Om$ is twice continuously differentiable. Let $0\le s \le
|\bo\Om|$ be an arc length parameter. Then the {\em curvature}
$\ka(s)$ is a continuous, nonnegative function on $\bo\Om$. We
will assume throughout the paper that $\Om$ is {\em strictly
convex in the sense of differential geometry}: $\ka>0$. In what
follows we refer to such $\Om$ as {\em regular billiard tables},
or              {\em regular convex domains}.

The elements of the {\em phase space of the billiard map} are the
inward pointing unit vectors $v$ based on $\bo\Om$. Let $0\le
\het\le\pi$ be the angle between $v$ and the positively oriented
$\bo\Om$. The coordinates $0\le s \le|\bo\Om|,0\le\het\le\pi$
induce a diffeomorphism of the phase space $Z$ and the cylinder
$\left(\R/|\bo\Om|\Z\right)\times[0,\pi]$.

\begin{figure}[htbp]
\begin{center}
\input{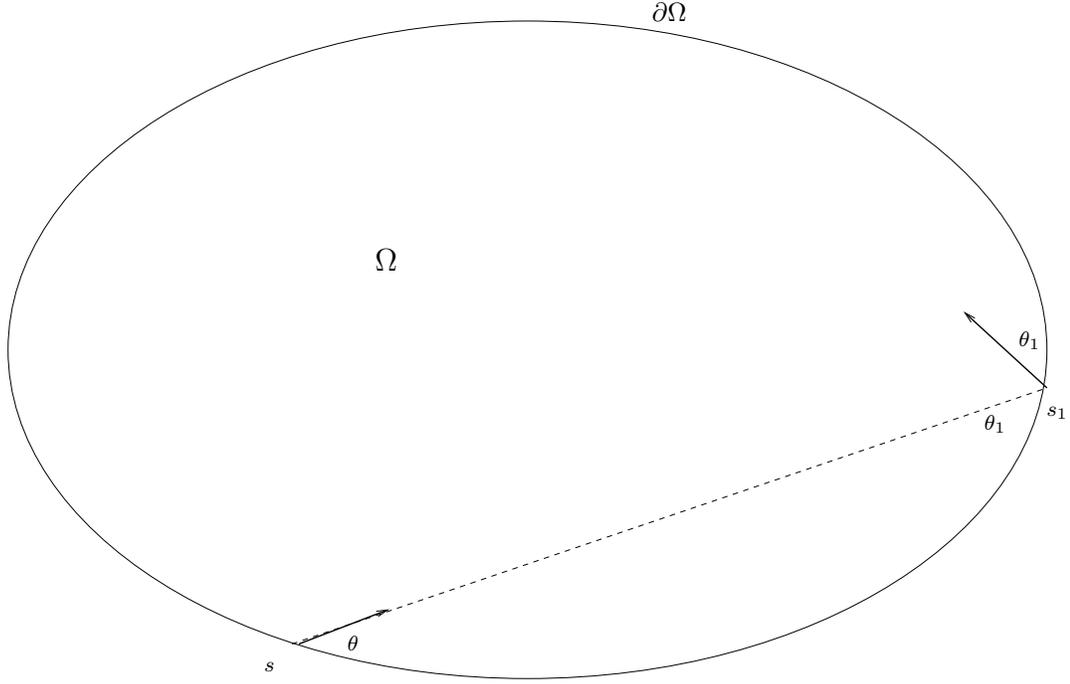}  
\caption{The billiard map for a regular convex domain.}
\label{bill_map_fig}
\end{center}
\end{figure}

The phase point $(s,\het)\in Z$ corresponds to the billiard ball
located at $s\in\bo\Om$, which is about to shoot of in the
direction that makes angle $\het$ with $\bo\Om$. This shot lands
at $s_1\in\bo\Om$. Let $0\le\het_1\le\pi$ be the other angle of
the chord $[s,s_1]$. The ball bounces elastically at the boundary
and is set to shoot of again. The {\em law of equal angles} yields
that the new vector $v_1$ makes angle $\het_1$ with $\bo\Om$. The
transformation $F:Z\to Z$ given by $F(s,\het)=(s_1,\het_1)$ is the
{\em billiard ball map} for $\Om$. Figure~~\ref{bill_map_fig}
illustrates the discussion.

In view of our assumptions on $\Om$, the billiard map is of class
$C^1$. Let $l(s,s_1)$ denote the length of the chord $[s,s_1]$.
The differential of the billiard map is given by the following
expressions \cite{GK}:
$$
\frac{\bo s_1}{\bo s}=\frac{\ka(s)l(s,s_1)-\sin\het}{\sin\het_1},\
\  \frac{\bo s_1}{\bo\het}=\frac{l(s,s_1)}{\sin\het_1},\
\frac{\bo\het_1}{\bo\het}=\frac{\ka(s_1)l(s,s_1)-\sin\het_1}{\sin\het_1},
$$
and
$$
\frac{\bo\het_1}{\bo
s}=\frac{\ka(s)\ka(s_1)l(s,s_1)-\ka(s)\sin\het_1-\ka(s_1)\sin\het}{\sin\het_1}.
$$

The billiard ball map is an {\em area preserving twist map}. The
classical results of Birkhoff on the dynamics of the billiard map
received a ``second life'' in the theory of area preserving twist
maps. See the accounts in \cite{Ba88}, \cite{MaFo94} and
\cite{KH}. We are concerned with a particular aspect of the
billiard ball map: The {\em invariant circles}.

\begin{defin}   \label{caustic_def}
Let $\Om$ be a regular billiard table. An {\em invariant circle}
for the billiard map on $\Om$ is a closed curve $\Ga\subset Z$
which is homotopic to a boundary component of $Z$ and is invariant
under the billiard map.
\end{defin}

By a theorem of Birkhoff, any invariant circle $\Ga$ is the graph
of a lipshitz function: $\het=h_{\Ga}(s)$. Thus, for every base
point $s\in\bo\Om$ there is a unique angle $\het=h_{\Ga}(s)$ such
that the ball shooting from $s$ in the direction $\het$ will
``stay'' on the invariant circle $\Ga$. For a typical $\Ga$ the
function $h_{\Ga}$ is not constant. See
\cite{Ba88,Gut94,GK,In88,GKn96,Laz73}. We will study invariant
circles such that $h_{\Ga}$ is constant. Both boundary components
of $Z$ are trivial invariant circles of that type. We will
consider only nontrivial invariant circles in what follows. To
simplify the terminology, we will often call them the {\em
invariant curves}. This is justified, since we will not study
other invariant curves.

\begin{figure}[htbp]
\input{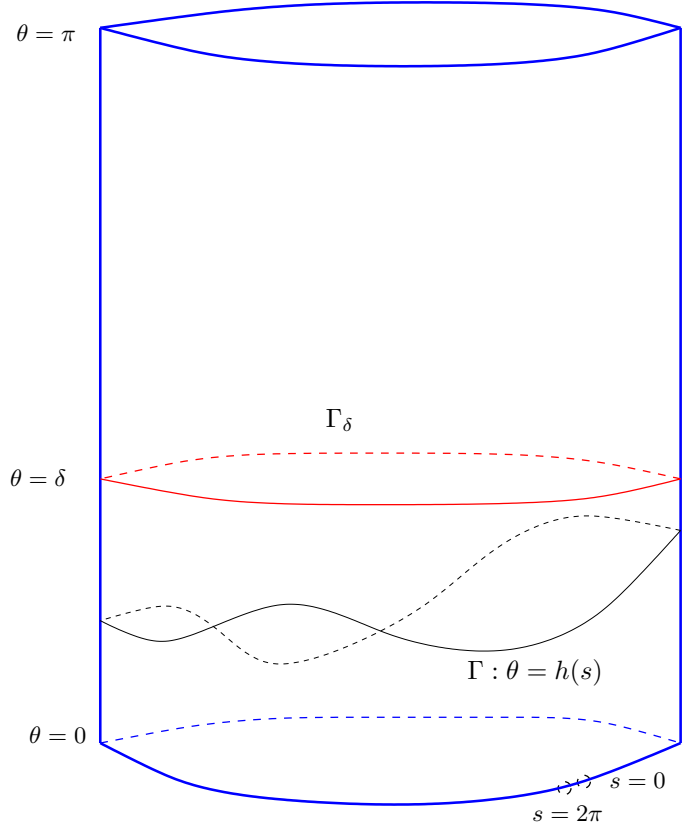} 
\begin{center}
\caption{Billiard map phase space with a general invariant circle
and a constant angle invariant circle.}
\label{inv_circ_fig}
\end{center}
\end{figure}

\begin{defin}   \label{const_type_def}
Let  $\Ga\subset Z$ be an invariant circle, and let
$\het=h_{\Ga}(s)$ the corresponding lipshitz function. If
$h_{\Ga}=\const$, we will say that $\Ga$ is a {\em constant angle
invariant circle}. A constant angle invariant circle is determined
by that angle, say $0<\de<\pi$. We will denote it by $\Ga_{\de}$.
See Figure~~\ref{inv_circ_fig}.
\end{defin}

\medskip

It is instructive to think of the phase space $Z$ as the space of
oriented lines (i. e., {\em rays}) intersecting $\Om$, or,
alternatively, as the space of directed chords in $\Om$. In this
representation, an invariant circle $\Ga$ is a one-parameter
family of rays. Its {\em envelope} $\ga\subset\RR$ is the {\em
caustic of $\Om$} corresponding to the invariant circle
$\Ga$.\footnote{The term ``caustic'' is widely used in the
geometric optics, mechanics, and the geometric theory of
singularities; in different contexts it stands for different,
although related things. We refer the reader to
\cite{Ar90,Be86,BG84} for many variations of this concept.}

Let $\Ga'$ be the family obtained from $\Ga$ by reversing the
directions of rays. Then $\Ga'$ is an invariant circle as well.
This is a consequence of the well known fact that the direction
reversing involution $\si:Z\to Z$ conjugates the billiard map with
its inverse: $\si F \si = F^{-1}$.\footnote{It is often called the
{\em billiard involution}.} It is clear that $\Ga$ and $\Ga'$ have
the same evolute; hence, the correspondence between invariant
circles and the caustics is $2$-to-$1$. The geometry of caustics
for regular billiard tables offers challenging open questions. See
\cite{Dou,GK} and \cite{GKn96} for this material. Since invariant
circles are determined by their caustics essentially uniquely, in
what follows we identify them; in particular, we will speak of
{\em general caustics} and of {\em constant angle caustics}.

\begin{rem}      \label{caustic_rem}
{\em The reader should keep in mind that the two invariant
circles, say $\Ga$ and $\Ga'=\si(\Ga)$ corresponding to the
caustic $\ga$ are distinct subsets of the phase space $Z$. Let
$0<r(\Ga)<1$ be the {\em rotation number} of the invariant circle.
Then $r(\Ga')=1-r(\Ga)$.

}
\end{rem}
%


Since $\bo\Om$ is strictly convex, we parameterize it by the
direction $0\le\al\le 2\pi$ of the tangent ray to $s\in\bo\Om$.
Thus, $s=s(\al)$. The derivative $\rho(\al)=ds/d\al$ is the {\em
radius of curvature} function for $\Om$. Set $\T=\R/2\pi\Z$. Then
the  billiard map is a diffeomorphism of $\T\times[0,\pi]$; we
will use the notation $F(\al,\het)=(\al_1,\het_1)$.

\begin{figure}[htbp]
\input{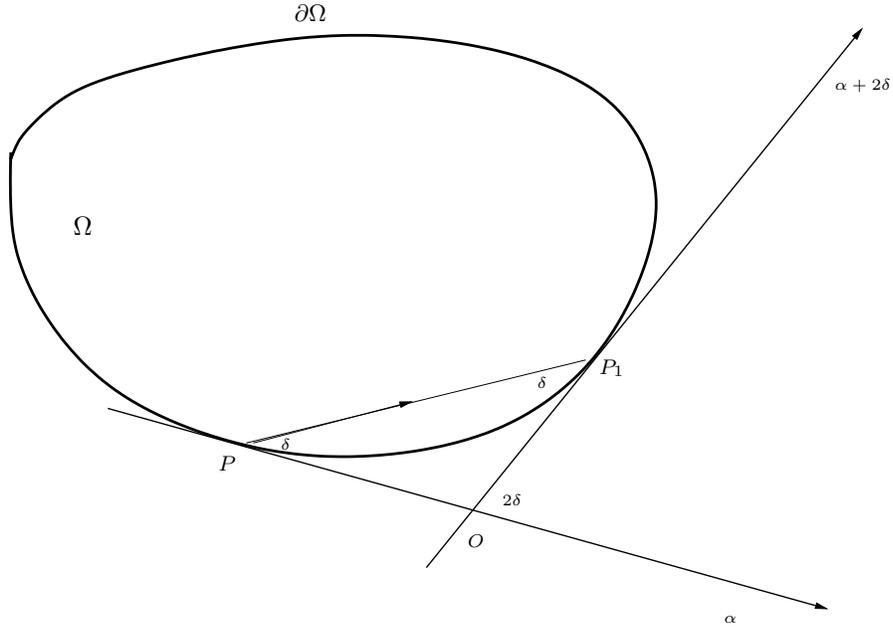} 
\begin{center}
\caption{The billiard map restricted to a constant angle invariant
circle.} \label{delta_fig}
\end{center}
\end{figure}

\begin{prop}             \label{const_caust_prop}
Let $\Om\subset\RR$ be a billiard table, and let
$\rho(\al),0\le\al\le 2\pi,$ be its radius of curvature. Then
$\Om$ has the constant angle caustic $\Ga_{\de}$ iff the function
$\rho(\cdot)$ satisfies the identity
\begin{equation}     \label{ident_eq1}
\int_{\al-\de}^{\al+\de}\rho(\xi)\sin(\al-\xi)d\xi=0.
\end{equation}
\begin{proof}
Set $P=P(\al)=(x(\al),y(\al))$ and let $P_1=P(\al_1)$. Let $O$ be
the intersection point of the tangent lines at $\al$ and $\al_1$.
From the triangle $POP_1$ we have $\al_1=\al+2\de$. See
Figure~~\ref{delta_fig}. As is well known
\begin{equation}     \label{rad_curv_eq}
x'(\al)=\rho(\al)\cos\al,\ y'(\al)=\rho(\al)\sin\al.
\end{equation}
Thus
$$
x(\al+2\de)-x(\al)=\int_{\al}^{\al+2\de}\rho(\xi)\cos\xi d\xi,
$$
$$
y(\al+2\de)-y(\al)=\int_{\al}^{\al+2\de}\rho(\xi)\sin\xi d\xi.
$$

The direction of the chord $[PP_1]$ is $\al+\de$. We introduce the
new variable $\be=\al+\de$. Thus, the slope of $[PP_1]$ is
$\tan\be$. Computing the slope from the coordinates of points $P$
and $P_1$, we obtain
\begin{equation}     \label{ident_eq2}
\frac{\int_{\be-\de}^{\be+\de}\rho(\xi)\sin\xi
d\xi}{\int_{\be-\de}^{\be+\de}\rho(\xi)\cos\xi d\xi}=\tan\be.
\end{equation}
Equation~~\eqref{ident_eq2} is an identity that holds for any
$\be\in\T$. Performing elementary trigonometric manipulations in
equation~~\eqref{ident_eq2}, and renaming the independent variable
by $\al$ again, we obtain the claim.
\end{proof}
\end{prop}

\medskip

\medskip

We will briefly review basic facts from harmonic analysis on the
circle. The reader may find proofs of the statements below in most
analysis textbooks.

If $g$ is a distribution on $\tf$, its Fourier transform is
defined by $\hg(n)=\int_{\tf}g(\al)e^{-in\al}d\al$ for   $n\in\Z$.
The radius of curvature has a Fourier expansion
\begin{equation}     \label{four_decom_eq}
\rho(\al)=\sum_{n\in\Z}c_ne^{in\al}
\end{equation}
where $c_n=\hat{\rho}(n)/2\pi$ are the Fourier coefficients. The
Fourier coefficients of a real function satisfy
$c_{-n}=\bar{c}_n$. Equation~~\eqref{four_decom_eq} is equivalent
to the trigonometric expansion
%
$$
\rho(\al)=a_0+\sum_{n\ge 1}a_n\cos n\al + b_n\sin n\al
$$
%
whose coefficients are real. The coefficients in these equations
are related by $a_0=c_0$ and $a_n=2\Re(c_n),b_n=-2\Im(c_n)$ for
$n>0$.

Denote by $x+y$ the group operation on $\tf$. Let $k(\cdot)$ be a
function or a distribution on $\tf$. The operator of convolution
with $k$ is defined by
\begin{equation}     \label{convol_eq}
(K\rho)(x)=\int_{\tf}\rho(x-\xi)k(\xi)d\xi=\int_{\tf}\rho(\xi)k(x-\xi)d\xi.
\end{equation}
The  standard notation for convolution operators is
$K(\rho)=\rho*k=k*\rho$.

Let $F_n$ be the complex line in the space of functions on $\tf$
spanned by $e^{inx}$. We view equation~~\eqref{four_decom_eq} as
the orthogonal decomposition by the subspaces $F_n$, $n\in\Z$.
Convolution operators preserve this decomposition. The restriction
$K|_{F_n}$ is the operator of multiplication by $\hat{k}(n)$. The
above discussion yields the following statement which is crucial
for Theorem~~\ref{const_caust_thm}.

\begin{lem}          \label{harm_claim}
Let $k(\cdot)$ be a distribution on $\tf$, and let
$\hat{k}(n),\,n\in\Z,$ be its Fourier transform. Let $K$ be the
operator of convolution  with the distribution $k(\cdot)$. Let
$\rho(\cdot)$ be a function on $\tf$, and let $c_n,\,n\in\Z,$ be
its Fourier coefficients.

\medskip

Then $K\rho=0$ iff $\hat{k}(n)c_n=0$ for all $n\in\Z$.
\end{lem}
\begin{thm}          \label{const_caust_thm}
Let $\Om\subset\RR$ be a regular, noncircular billiard table. Let
$\rho(\cdot)$ be the radius of curvature of $\bo\Om$, and let
$c_n,\,n=1,2\dots$ be its Fourier coefficients. Then $\Om$ has the
constant angle caustic $\Ga_{\de}$  iff the following conditions
hold:

\begin{itemize}
\item[i)] There exist  $n>1$ such that
\begin{equation}     \label{sin_eq}
\frac{\sin(n-1)\de}{n-1}=\frac{\sin(n+1)\de}{n+1};
\end{equation}
\item[ii)] We have $c_k=0$ for all $k>1$ such that
equation~~\eqref{sin_eq} is not satisfied.
\item[iii)] We have $c_n\ne0$ for at least one $n>1$ such that
equation~~\eqref{sin_eq} is satisfied.
\end{itemize}
\begin{proof}
By Proposition~~\ref{const_caust_prop}, $\Ga_{\de}$ is a caustic
for $\Om$ iff $\rho(\cdot)$ belongs to the zero space of the
convolution with the function
\begin{equation}     \label{kernel_eq}
k(x)=(\sin x)\,1_{[-\de,\de]}.
\end{equation}
The function $k(\cdot)$ is odd, hence $\hat{k}(0)=0$. By a
straightforward computation, for $n>1$ we have
$$
i\cdot
\hat{k}(n)=\frac{\sin(n-1)\de}{n-1}-\frac{\sin(n+1)\de}{n+1}.
$$

\medskip

It is well known that for any billiard table $\Om$ we have
$c_1=0$. By Lemma~~\ref{harm_claim}, $\Ga_{\de}$ is a caustic for
$\Om$ iff $c_m\hat{k}(m)=0$ for all $m$. By the preceding
discussion, $\Ga_{\de}$ is a caustic iff
$$
c_m\left[\frac{\sin(m-1)\de}{m-1}-\frac{\sin(m+1)\de}{m+1}\right]=0
$$
for $m>1$. Since $\Om$ is not circular, at least one coefficient,
say $c_n$, does not vanish. But $\Ga_{\de}$ being a caustic,
equation~~\eqref{sin_eq} holds for all $n>1$ such that $c_n\ne0$.
\end{proof}
\end{thm}

\section{Constant angle caustics and a chain of trigonometric equations}    \label{reduce}
By Theorem~~\ref{const_caust_thm}, the description of billiard
tables with constant angle caustics hinges on solving
equation~~\eqref{sin_eq}.  In this section we will reduce
equation~~\eqref{sin_eq} to a chain of trigonometric equations
involving the function $\tan(\cdot)$.

Recall that $\de\in A$ if there exists $n>1$ such that
equation~~\eqref{sin_eq} holds. Let $A_n\subset A$ be the set of
$\de\in(0,\pi)$ such that the pair $\de,n$ satisfies
equation~~\eqref{sin_eq}. Thus
$$
A=\cup_{n=2}^{\infty}A_n.
$$
\begin{lem}             \label{trig_eqn_lem}
Let $n>1$. Then the following claims hold.
\begin{itemize}
\item[i)] We have $\frac{\pi}{2}\in A_n$ iff $n$ is odd.
\item[ii)] Set $\tA_n=A_n\setminus\{\frac{\pi}{2}\}$. Then $\tA_n$
is the set of solutions in $(0,\pi)$ of the equation $\tan
n\de=n\tan\de$.
\end{itemize}
\begin{proof}
Set $\de=\frac{\pi}{2}$. If $n$ is odd, then both sides in
equation~~\eqref{sin_eq} vanish, hence $\frac{\pi}{2}\in A_n$. If
$n$ is even, then the numerators in equation~~\eqref{sin_eq} are
$\pm1$, and their signs are opposite. Thus, $\frac{\pi}{2}\notin
A_n$, proving claim i).

Let $\de\in A_n$. Arguing as above, we establish that
$\sin(n+1)\de=0$ iff $n$ is odd and $\de=\frac{\pi}{2}$. Hence for
$\de\in\tA_n$ we have $\sin(n-1)\de,\sin(n+1)\de\ne 0$. Therefore,
$\tA_n$ is the set of $\de\in(0,\pi)$ satisfying
$$
\frac{\sin(n-1)\de}{\sin(n+1)\de}=\frac{n-1}{n+1}.
$$
We rewrite this as
\begin{equation}     \label{trig_eq}
\frac{\sin n\de\cos\de-\cos n\de\sin\de}{\sin n\de\cos\de+\cos
n\de\sin\de}=\frac{n-1}{n+1}.
\end{equation}
If $\cos n\de=0$, then the left hand side in
equation~~\eqref{trig_eq} is $1$, which is impossible. Thus, $\cos
n\de\ne 0$. Dividing the numerators and the denominators in
equation~~\eqref{trig_eq} by $\cos\de\cos n\de$, we obtain
$$
\frac{\tan n\de - \tan\de}{\tan n\de + \tan\de}=\frac{n-1}{n+1}.
$$
Claim 2 follows.
\end{proof}
\end{lem}

\medskip

For $X\subset\R$ and $a\in\R$ let $\{a-X\}=\{a-x:x\in X\}$. Set
$\tA=A\setminus\{\frac{\pi}{2}\}$. Then
\begin{equation}   \label{const_angles_eq}
\tA=\cup_{n=2}^{\infty}\tA_n.
\end{equation}
Set $B_n=A_n\cap(0,\pi/2)$ and $B=A\cap(0,\pi/2)$.
Lemma~~\ref{trig_eqn_lem} and the preceding discussion imply the
following.
\begin{prop}          \label{angles_prop}
Let $n>1$. Then for $n$ even, $A_n=B_n\cup\{\pi-B_n\}$ and for $n$
odd, $A_n=B_n\cup\{\pi-B_n\}\cup\{\frac{\pi}{2}\}$. Moreover,
$B_n$ is the set of solutions in $(0,\pi/2)$ of the equation
\begin{equation}   \label{tan_nx_eq}
\tan n x=n\tan x.
\end{equation}
\end{prop}

\section{Analysis of trigonometric equations}    \label{trig_qual}
\begin{figure}[htbp]
\input{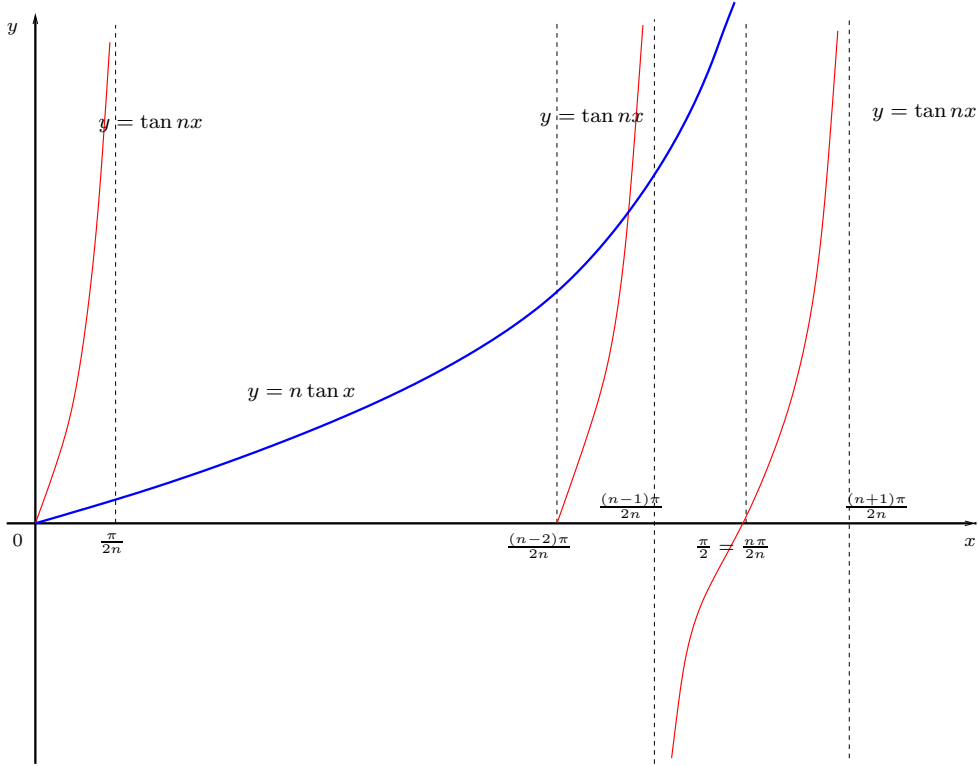} 
\begin{center}
\caption{The graphs of functions $y=\tan nx$ and $y=n\tan x$ for
$n$ even.} \label{even_n_fig}
\end{center}
\end{figure}

\begin{figure}[htbp]
\input{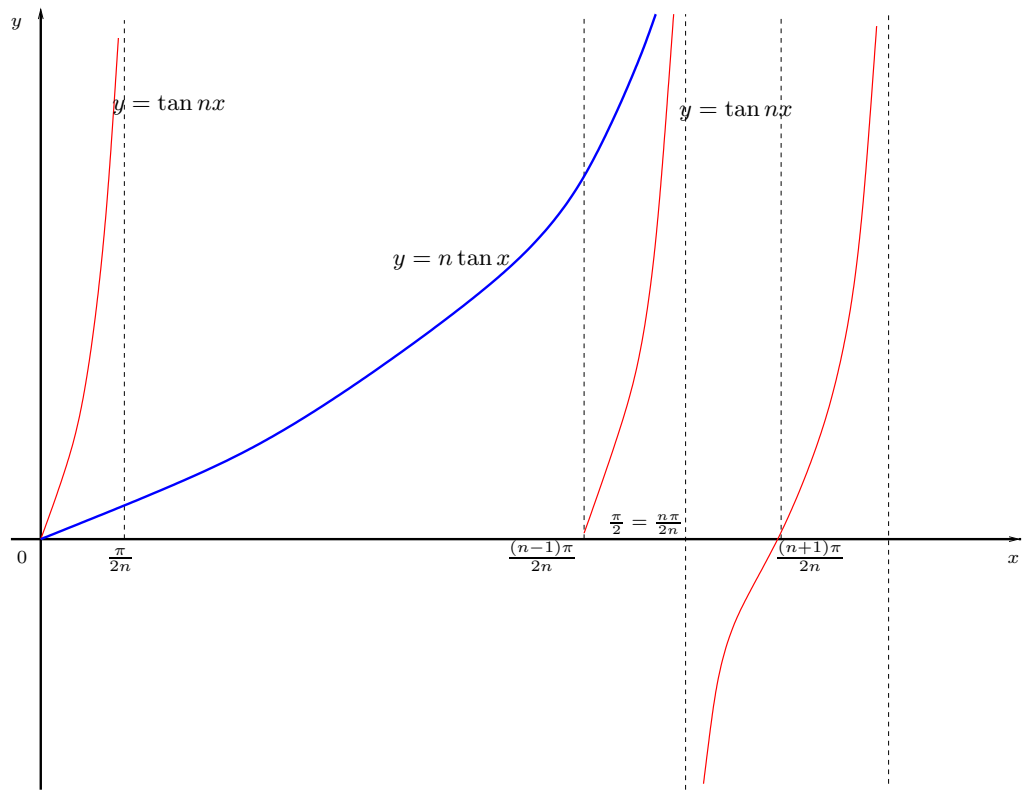} 
\begin{center}
\caption{The graphs of functions $y=\tan nx$ and $y=n\tan x$ for
$n$ odd.} \label{odd_n_fig}
\end{center}
\end{figure}

In this section we begin to analyze solutions of the chain of
equations~~\eqref{tan_nx_eq} in the interval $(0,\pi/2)$.
\begin{prop}         \label{qual_prop}
{\em 1}. Let $n>1$ be even. Then $B_n$ consists of $\frac{n}{2}-1$
points $\xi_k^{(n)}$, where
$$
\frac{2k}{2n}\pi<\xi_k^{(n)}<\frac{(2k+1)}{2n}\pi:\,k=1,\dots,\frac{n}{2}-1.
$$

\noindent {\em 2}. Let $n>1$ be odd. Then $B_n$ consists of
$\frac{n-1}{2}-1$ points $\xi_k^{(n)}$, where
$$
\frac{2k}{2n}\pi<\xi_k^{(n)}<\frac{(2k+1)}{2n}\pi:\,k=1,\dots,\frac{n-1}{2}-1.
$$
\begin{proof}
The graph of the function $y=\tan nx$ on $(0,\pi/2)$ is the
disjoint union of $n$ connected curves; we will call them {\em
branches}. A branch is defined on the interval
$\frac{k}{2n}\pi<x<\frac{k+1}{2n}\pi:\,0\le k \le n-1$. Set
$I_k^{(n)}=(\frac{k}{2n}\pi,\frac{k+1}{2n}\pi)$. Each branch
extends by continuity to one of the endpoints of $I_k^{(n)}$.
These endpoints don't enter in our analysis, and we ignore them in
what follows. We say that a branch is {\em positive} (resp. {\em
negative}) if it belongs the the upper (resp. lower) halfplane.

Positive branches correspond to $I_k^{(n)}$ with $k$ even. Thus,
there are $n/2$ (resp. $(n-1)/2$) positive branches if $n$ is even
(resp. odd). We observe that each point in $B_n$ belongs to the
intersection of the graph of $y=n\tan x$ on $(0,\pi/2)$ with a
positive branch; this intersection contains at most one point. See
Figure~~\ref{even_n_fig} and Figure~~\ref{odd_n_fig}.

Comparing the asymptotics of $n\tan x$ and $\tan nx$ as $x\to 0+$,
we  see that the first branch, which  corresponds to $k=0$,  does
not yield an intersection point. When $n$ is even, all other
positive branches intersect the graph $n\tan x$. This proves claim
1. Let now $n$ be odd. Then both the last branch and the graph of
$y=n\tan x$ are asymptotic to the vertical line $x=\pi/2$.
Comparing the asymptotics of $n\tan x$ and $\tan nx$ as
$x\to\frac{\pi}{2}-$, we see that the curves do not intersect.
This proves claim 2. We leave details to the reader.
\end{proof}
\end{prop}
%


%

We will state an immediate consequence of
Proposition~~\ref{qual_prop}. Recall that $A_n\subset(0,\pi)$ is
the set of numbers satisfying equation~~\eqref{sin_eq} and that
$A=\cup_{n>1}A_n$.
\begin{corol}               \label{qual_cor}
We have $|A_n|=n-2$. The sets $A_n$ are $2\pi/n$ dense in
$(0,\pi)$.
\begin{proof}
By Lemma~~\ref{trig_eqn_lem} and Proposition~~\ref{angles_prop},
$|A_n|=2|B_n|$ if $n$ is even and $|A_n|=2|B_n|+1$ if $n$ is odd.
The first claim now follows from Proposition~~\ref{qual_prop}. The
above propositions imply that the distances between consecutive
points of $A_n$ are at most $\pi/n$. Besides, the distances from
$A_n$ and the endpoints of $(0,\pi)$ are at most $2\pi/n$.
\end{proof}
\end{corol}

\section{Immediate implications}
\label{impli}
The results of section~~\ref{birkhoff}, section~~\ref{reduce}, and
section~~\ref{trig_qual} have immediate consequences for the
billiard  and for the floating. We begin with the former. We will
say that the billiard tables $\Om_1,\Om_2$ are {\em conformally
equivalent} if there is a conformal mapping $F:\RR\to\RR$ such
that $F(\Om_1)=\Om_2$. For instance, all discs in $\RR$ are
conformally equivalent.

\begin{corol}      \label{1_param_prop}
There is a dense countable set
$\tA\subset(0,\pi)\setminus\{\pi/2\}$ such that the following
holds.

\noindent {\em 1}. For any $\de\in\tA$ there is $n>1$ and a real
analytic $1$-parameter family $\Om_{n,\tau},\,0\le\tau<1,$ of
conformally inequivalent, regular billiard tables having the
constant angle caustic $\Ga_{\de}$. The curves $\bo\Om_{n,\tau}$
are real analytic; $\bo\Om_{n,0}$ is the unit circle.

\noindent {\em 2}. A regular billiard table $\Om$ has the caustic
$\Ga_{\pi/2}$ iff $\bo\Om$ is a curve of constant width.

\noindent {\em 3}. Let $0<\de<\pi$ belong to the complement of
$\tA\cup\{\pi/2\}$ in $(0,\pi)$. If a regular billiard table $\Om$
has the constant angle caustic $\Ga_{\de}$, then $\Om$ is
circular.
\begin{proof}
Let $\tA_n$ and $\tA$ be as in equation~~\eqref{const_angles_eq}.
Then $\de\in\tA$ iff there exists $n>1$ such that $\tan n
\de=n\tan\de$. Let $a,b\in\R$ be arbitrary. Set
$$
\rho(\al)=1+a\cos n\al + b\sin n\al.
$$
By elementary trigonometry, there exists $\al_0$ depending on
$a,b,n$ such that $\rho(\al)=1+\sqrt{a^2+b^2}\sin n(\al+\al_0)$.
This is the radius of curvature of a regular billiard table iff
$a^2+b^2<1$. Different values of $\al_0$ correspond to isometric
billiard tables. Set $\rho_{n,\tau}(\al)=1+\tau\sin n\al$.


Integrating
equation~~\eqref{rad_curv_eq}, we obtain
\begin{equation}        \label{table_eq1}
x_{n,\tau}(\al)=\xi_0+\sin\al+\frac{\tau}{2(n-1)}\cos(n-1)\al-\frac{\tau}{2(n+1)}\cos(n+1)\al,
\end{equation}
\begin{equation}         \label{table_eq2}
y_{n,\tau}(\al)=\eta_0-\cos\al+\frac{\tau}{2(n-1)}\sin(n-1)\al-\frac{\tau}{2(n+1)}\sin(n+1)\al.
\end{equation}
Fixing the constants $\xi_0,\eta_0$, we obtain a real analytic
family $\Om_{n,\tau}$. This proves claim 1. Claim 2 will follow
from the discussion in the beginning of
section~~\ref{const_width}. Claim 3 is immediate from
Theorem~~\ref{const_caust_thm}.
\end{proof}
\end{corol}

Set $\rho(\al)=c+a\cos n\al + b\sin n\al$. This formula, provided
$0\le\sqrt{a^2+b^2}<c$, yields a $3$-parameter family of functions
that serve as radii of curvature for billiard tables $\Om$ having
the caustic $\Ga_{\de}$. By
equations~~\eqref{table_eq1},~~\eqref{table_eq2}, we have a
$5$-parameter family of these domains. However, the conformal
equivalence eats up $4$ of the parameters.
We can view equations~~\eqref{table_eq1},~~\eqref{table_eq2} as a
deformation $\Om_{n,\tau}$ of the circular table.

\medskip

The following is the counterpart of Corollary~~\ref{1_param_prop}
for the floating in neutral equilibrium. Its claims are the
reformulations of the corresponding claims in
Corollary~~\ref{1_param_prop}; we do not repeat the proof.
\begin{corol}      \label{float_cor}
There is a dense countable set
$\tA\subset(0,\pi)\setminus\{\pi/2\}$ such that the following
holds.

\noindent {\em 1}. For any $\de\in\tA$ there is $n>1$ and a real
analytic $1$-parameter family $\Om_{n,\tau},\,0<\tau<1,$ of
conformally inequivalent planar domains with real analytic
boundaries that float in neutral equilibrium at any orientation
with the contact angle $\pi-\de$. The domain $\Om_{n,0}$ is the
unit circle.

\noindent {\em 2}. A regular convex domain floats in neutral
equilibrium at any orientation with the contact angle $\pi/2$ if
and only if its boundary is a curve of constant width.

\noindent {\em 3}. Let $0<\de<\pi$ belong to the complement of
$\tA\cup\{\pi/2\}$ in $(0,\pi)$. If a regular  convex domain
floats in neutral equilibrium at any orientation with the contact
angle $\de$, then it is a disc.
\end{corol}


Recall that $B_n\subset(0,\pi/2)$ is the set of solutions to
equation~~\eqref{tan_nx_eq}.  To continue our study of floating in
neutral equilibrium, we need further number theoretic information
about these sets. Below we formulate questions about floating
and/or billiard whose answers depend on this information.


%
\begin{quest}   \label{number_quest1}
{\em Let $\de\in B$. Describe the set of billiard tables $\Om$
such that $\Ga_{\de}$ is a caustic. Equivalently, describe the set
of cross-sections of cylinders that float in neutral equilibrium
with the contact angle $\pi-\de$ at any orientation.

}
\end{quest}
%


%
\begin{quest}   \label{number_quest2}
{\em Let $\de\in B$. Let $\Om$ be a billiard table with the
caustic $\Ga_{\de}$. Let $\vp:\Ga_{\de}\to\Ga_{\de}$ be the
restriction of the billiard map on $\Om$ to $\Ga_{\de}$. Can $\vp$
be periodic?

}
\end{quest}

In order to answer Question~~\ref{number_quest1}, we need to
investigate the intersections $B_n\cap B_m$ for $m\ne n$. In
particular, we need to know for what pairs $m\ne n$ the set
$B_n\cap B_m$ is nonempty. The answer to Q
uestion~~\ref{number_quest2} depends on whether $\de\in\pi\QQ$ or
not. In particular, if $B\cap\pi\QQ$ is empty, the answer is
negative. In section~~\ref{trig_quant} and the following sections
we will study the sets of solutions to equation~~\eqref{tan_nx_eq}
and equation~~\eqref{sin_eq}. The solution $\de=\pi/2$ of the
latter is special. The corresponding planar domains have been
studied by geometers from an independent viewpoint. We briefly
review this in the next section.

\section{The   caustics $\Ga_{\pi/2}$ and curves of constant width}    \label{const_width}
To illustrate the preceding discussion, we will now study the
question: Which billiard tables have the caustic $\Ga_{\pi/2}$?
Let $\Om$ be a regular billiard table. Then $\Ga_{\pi/2}$ is a
caustic iff any chord which is perpendicular to $\bo\Om$ at one of
its ends, is also perpendicular to $\bo\Om$ at the other end. The
values of the angle parameter at these points are $\al,\al+\pi$.
The length $d(\al)=|P(\al)P(\al+\pi)|$ is the {\em width of $\Om$
in the direction $\al+\pi/2$}. The chord $[P(\al)P(\al+\pi)]$ is
perpendicular to $\bo\Om$ iff  $\al+\pi/2$ is a critical point for
the function $d(\cdot)$ . Therefore, $\Ga_{\pi/2}$ is a caustic
iff $d(\al)=\const$. These curves are known in geometry as the
{\em curves of constant width} \cite{BY60}. Thus, a regular
billiard table $\Om$ has the caustic $\Ga_{\pi/2}$ iff $\Om$ is a
domain of constant width.

\begin{figure}[htbp]
\input{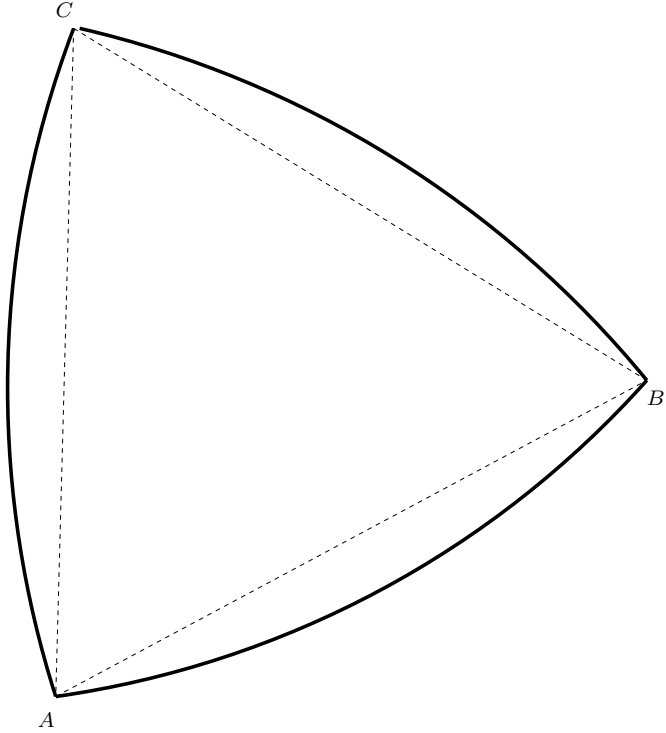} 
\begin{center}
\caption{Reuleaux triangle: A domain of constant width with
corners.} \label{releux_fig}
\end{center}
\end{figure}

We point out that the analysis below assumes that $\bo\Om$ is
twice continuously differentiable. In particular, it is not valid
for domains of constant width with corners. The boundary of the
famous example of such a domain, the Reuleaux triangle
\cite{BY60}, consists of three circular arcs of the same radius;
it has corners at the endpoints of the arcs. See
Figure~~\ref{releux_fig}. The Reuleaux triangle is not a regular
billiard table.

\begin{corol}             \label{const_width_cor}
Let $\Om$ be a regular billiard table, and let $\rho(\cdot)$ be
its radius of curvature. Then $\Ga_{\pi/2}$ is a caustic for $\Om$
iff we have the identity
\begin{equation}     \label{const_eq}
\rho(\al)+\rho(\al+\pi)=\const.
\end{equation}
\begin{proof}
Let $c_m,m\in\Z,$ be the Fourier coefficients of $\rho$. By the
proof of Theorem~~\ref{const_caust_thm}, $\Ga_{\pi/2}$ is a
caustic iff
$$
c_m\left[\frac{\sin\frac{(m-1)\pi}{2}}{m-1}-\frac{\sin\frac{(m+1)\pi}{2}}{m+1}\right]=0
$$
for all $m>1$. For $m=2k$ this means $4kc_{2k}/(4k^2-1)=0$,
yielding $c_{2k}=0$. For odd $m$ the equation holds for any $c_m$.
Thus, the caustic $\Ga_{\pi/2}$ exists iff the radius of curvature
has the Fourier expansion of the form
\begin{equation}     \label{odd_eq}
\rho(\al)= c_0+\sum_{m\,\mbox{odd}}c_me^{im\al}.
\end{equation}
Set $\rho_0(\al)=\rho(\al)-c_0$. Then equation~~\ref{odd_eq} holds
iff $\rho_0$ is an odd function on $\tf$. Equivalently,
$\rho(\al+\pi)+\rho(\al)=2c_0$.
\end{proof}
\end{corol}

\begin{rem}     \label{const_width_rem}
{\em We point out that the identity equation~~\eqref{odd_eq}
characterizes all billiard tables $\Om$ with the caustic
$\Ga_{\pi/2}$, including  the circular billiard table. By the
discussion preceding Corollary~~\ref{const_width_cor}, the width
of any such $\Om$ is constant, and is equal to $2c_0$. Let
$|\bo\Om|$ be the perimeter of $\Om$. If $\Om$ has constant width,
we denote it by $w(\Om)$. By the above argument, for a curve of
constant width we have
$$
\rho(\al+\pi)+\rho(\al)=w(\Om).
$$
Integrating this equation and using that
$\int_{\tf}\rho(\al)d\al=|\bo\Om|$, we obtain the identity
\begin{equation}     \label{width_eq}
\pi\cdot w(\Om)=|\bo\Om|.
\end{equation}

Note that we have used the regularity of $\bo\Om$ to derive
equation~~\eqref{width_eq}. In fact, it is valid for arbitrary
curves of constant width; it is called Barbier's theorem. Another
amusing fact about domains of  constant width is the
Blaschke-Lebesgue theorem \cite{BY60}. It says that amongst the
domains of a fixed constant width the Reuleaux triangle has the
smallest area. By the isoperimetric theorem, the disc has the
biggest area. Let $\Om$ be any domain of constant width $w$; let
$|\Om|$ be the area of $\Om$. By an elementary calculation
$$
\frac{\pi-\sqrt{3}}{2}w^2\le|\Om|\le\frac{\pi}{4}w^2.
$$
The equalities take place only for the Reuleaux triangle and the
disc.

}
\end{rem}

\section{Trigonometric equations and a family of polynomials}    \label{trig_quant}
We will now obtain quantitative information about the solutions of
equations~~\eqref{tan_nx_eq}.
\begin{lem}     \label{prelim_lem}
Let $n\ge 1$. There are polynomials $P_n,Q_n$ such that
\begin{equation}   \label{rat_func_eq}
\tan nx =\frac{P_n(\tan x)}{Q_n(\tan x)}.
\end{equation}
Polynomials $P_n,Q_n$ are uniquely determined by the recurrence
relations
\begin{equation}   \label{recurr_eq}
P_{n+1}(z)=P_n(z)+zQ_n(z),\ Q_{n+1}(z)=Q_n(z)-zP_n(z)
\end{equation}
and the initial data $P_1(z)=z,\ Q_1(z)=1$. The polynomial $P_n$
(resp. $Q_n$) is odd (resp. even). The degree of each of the two
polynomials is either  $n$ or $n-1$, depending on the parity of
$n$.
\begin{proof}
The formula
$$
\tan(x+y)=\frac{\tan x + \tan y}{1-\tan x\tan y}
$$
in the special case $y=nx$ yields
$$
\tan(n+1)x=\frac{\tan nx + \tan x}{1-\tan nx\tan x}.
$$
The claims follow by induction on $n$.
\end{proof}
\end{lem}
\begin{rem}     \label{chebyshev_rem}
{\em Polynomials $P_n,Q_n$ can be expressed in terms of the
Chebyshev polynomials of the first and the second kind. We will
not pursue this approach here.

}
\end{rem}

\begin{prop}    \label{polynom_prop}
The polynomials in equation~~\eqref{rat_func_eq} satisfy
\begin{equation}   \label{P_n_eq}
-2P_n(z)=i^{n+1}(z-i)^n+(-i)^{n+1}(z+i)^n
\end{equation}
and
\begin{equation}   \label{Q_n_eq}
2Q_n(z)=i^n(z-i)^n+(-i)^n(z+i)^n.
\end{equation}
\begin{proof}
We rewrite equation~~\eqref{recurr_eq} as
%
$$
\begin{bmatrix}P_{n+1}(z)\\Q_{n+1}(z)\end{bmatrix}=
\begin{bmatrix}1&z\\-z&1\end{bmatrix}
\begin{bmatrix}P_n(z)\\Q_n(z)\end{bmatrix}.
$$
%
The claims follow by the elementary algebra.
\end{proof}
\end{prop}

\begin{corol}     \label{trig_polynom_cor}
For $n>1$ set
\begin{equation}   \label{R_n_eq}
R_n(z)=-\frac12\,i^n\left[(nz+i)(z-i)^n+(-1)^n(nz-i)(z+i)^n\right].
\end{equation}

Let $0<x<\pi/2$, and set $z=\tan x$. Then $x\in B_n$ iff $z$ is a
positive root of the polynomial $R_n$.
\begin{proof}
By Proposition~~\ref{angles_prop}, equation~~\eqref{tan_nx_eq},
and Lemma~~\ref{prelim_lem}, $x\in B_n$ iff $\tan x=z>0$ satisfies
$ P_n(z)-nzQ_n(z)=0$. By equations~~\eqref{P_n_eq}
and~~\eqref{Q_n_eq}, $P_n(z)-nzQ_n(z)=R_n(z)$.
\end{proof}
\end{corol}
%


\section{Number theoretic conjectures and implications}    \label{numb_theor}
We have reduced our investigation of equations~~\eqref{tan_nx_eq}
to a study of roots of the polynomials $R_n$. We now continue to
study these polynomials, and bring in some number theory.

\subsection{Polynomials and fractional linear transformations}
\label{frac_lin_sub}
\hfill \break  We will investigate the roots of $R_n$. The
following lemma summarizes the immediate properties of these
polynomials.

\begin{lem}     \label{elementary_lem}
Let $n\ge 1$. Then the following holds:
\begin{itemize}
\item[i)] The polynomials $R_n$ are real, odd polynomials;
\item[ii)] The degree of  $R_n$ is equal to $n+1$ for $n$ even, and to $n$ for $n$ odd;
\item[iii)] The highest coefficient of $R_n$ is $2n$ for $n$ even and $\pm1$ for $n$ odd;
\item[iv)] We have $R_n(z)=O(z^3)$;
\item[v)] The roots of $R_n$ are real and simple, except for the zero root, which has multiplicity three.
\end{itemize}
\begin{proof}
Claims i) - iv) follow either from $R_n(z)=P_n(z)-nzQ_n(z)$ or
directly from equation~~\eqref{R_n_eq}. We will prove claim v).
Suppose $n$ is even; set $n=2k$. Then $\deg(R_n)=2k+1$. By iv),
$R_n$ has at most $2k-2$ nonzero roots, counted with
multiplicities. By claim 1 in Proposition~~\ref{qual_prop} and
Corollary~~\ref{trig_polynom_cor}, $R_n$ has $k-1$ distinct
positive roots. By i), $R_n$ has $k-1$ distinct negative roots,
hence the claim. The case of odd $n$ is similar, and we leave it
to the reader.
\end{proof}
\end{lem}

Let $\begin{bmatrix}a&b\\c&d\end{bmatrix}$ be a nondegenerate
matrix. We will use the notation
%
$$
\begin{bmatrix}a&b\\c&d\end{bmatrix}\circ z=\frac{az+b}{cz+d}.
$$
%
Let $A,A'$ be nondegenerate matrices. We will write $A\sim A'$ to
mean that $A'A^{-1}$ is a scalar matrix. Then $A\sim A'$ holds iff
$A\circ z\equiv A'\circ z$.
\begin{prop}     \label{roots_prop}
Let $n>1$. There is a $1$-to-$1$ correspondence, preserving the
multiplicities, between the nonzero roots of $R_n$ and the roots
of the equation
\begin{equation}     \label{funny_eq}
\zeta^n=(-1)^{n+1}\begin{bmatrix}n+1&n-1\\n-1&n+1\end{bmatrix}\circ
\zeta,
\end{equation}
other than $\zeta=\pm 1$.
\begin{proof}
By equation~~\eqref{R_n_eq}, we have $R_n(z)=0$ iff
\begin{equation}     \label{before_funny_eq}
\left(\frac{z-i}{z+i}\right)^n=\left(-1\right)^{n+1}\frac{nz-i}{nz+i}.
\end{equation}

We recall a few well known facts. The fractional linear
transformations
$$
\zeta=\begin{bmatrix}1&-i\\1&i\end{bmatrix}\circ z,\
z=\begin{bmatrix}i&i\\-1&1\end{bmatrix}\circ\zeta.
$$
are inverse to each other; they induce a diffeomorphism of
$\R\cup\infty$ onto the unit circle which sends the natural
orientation of the real axis to the counter clockwise orientation
of the unit circle.

Setting $F_n(z)=\frac{nz-i}{nz+i}$, we rewrite equation~~
\eqref{before_funny_eq} as
$$
\left(F_1(z)\right)^n=\left(-1\right)^{n+1}F_n(z).
$$

Setting $F_1(z)=\zeta, z=F_1^{-1}(\zeta)$, and using that
$$
\begin{bmatrix}n&-i\\n&i\end{bmatrix}
\begin{bmatrix}i&i\\-1&1\end{bmatrix}=\begin{bmatrix}ni+i&ni-i\\ni-i&ni+i\end{bmatrix}\sim\begin{bmatrix}n+1&n-1\\n-1&n+1\end{bmatrix},
$$
we obtain equation~~\eqref{funny_eq}.

We have proved that the transformation
$\zeta=\begin{bmatrix}1&-i\\1&i\end{bmatrix}\circ z$ induces a
multiplicity preserving isomorphism between the roots $z$ of $R_n$
such that $F_1(z)\ne\infty$ and the solutions $\zeta\ne
F_1(\infty)$ of equation~~\eqref{funny_eq}. Using that
$F_1(0)=-1,F_1(\infty)=1$, and the information about the roots of
$R_n$ contained in Lemma~~\ref{elementary_lem}, we obtain the
claim.
\end{proof}
\end{prop}
\begin{rem}     \label{minus_one_rem}
{\em Proposition~~\ref{roots_prop} singles out the roots
$\zeta=\pm1$ of equation~~\eqref{funny_eq}. Observe that $-1$ is
always a root of multiplicity three for this equation, while $1$
is a (simple) root iff $n$ is odd. To explain this, we note that
$-1=F_1(0)$, while $1=F_1(\infty)$. Observe that $0$ is a
multiplicity three root of $R_n$; the appearance of $\infty$ as a
``root'' of $R_n$ is due to the circumstance that in the beginning
of the proof of Proposition~~\ref{roots_prop} we have put the
equation $R_n(z)=0$ in the form
\begin{equation}   \label{again_R_n_eq}
(nz+i)(z-i)^n=(-1)^{n+1}(nz-i)(z+i)^n.
\end{equation}
If $n$ is odd, the leading terms in both sides of
equation~~\eqref{again_R_n_eq} have the same sign when
$z\to\infty$; if $n$ is even, the signs are opposite.

}
\end{rem}

Equation~~\eqref{funny_eq} involves a rational function whose
denominator is $(n-1)\zeta+(n+1)$. Getting rid of the denominator
and using the variable $x=-\zeta$, we obtain an equivalent
polynomial equation:
$$
-(n-1)x^{n+1}+(n+1)x^n-(n+1)x+(n-1)=0.
$$
The two corollaries below follow immediately from
Proposition~~\ref{roots_prop} and the preceding discussion.

\begin{corol}      \label{filaseta_cor}
Let $n>1$. Set
\begin{equation}   \label{filaseta_eq}
S_n(x)=(n-1)\left[x^{n+1}-1\right]-(n+1)\left[x^n-x\right].
\end{equation}
Then all roots of the polynomials $S_n$ belong to the unit circle
$\{|x|=1\}$. The number $1$ is a root of multiplicity three. The
number $-1$ is a simple root of $S_n$ if $n$ is odd, and
$S_n(-1)\ne 0$ if $n$ is even. The remaining roots of $R_n$ are
simple.
\end{corol}

In what follows we will refer to the roots $x\ne\pm1$ of $S_n$ as
the {\em complex roots}.

\begin{corol}      \label{gutkin_cor}
Let $n>1$. The transformation $z\mapsto x$ given by
$$
x=-\frac{z-i}{z+i}
$$
induces a $1$-to-$1$ correspondence between the nonzero  roots of
the polynomial $R_n$ and the complex roots of the polynomial
$S_n$. Moreover, this transformation sends the positive (resp.
negative) roots of $R_n$ to the roots of $S_n$ such that $\Im x
>0$ (resp. $\Im x < 0$).
\end{corol}

\subsection{Polynomials $S_n$: Conjectures and supporting evidence}
\label{conj_sub}
\hfill \break Let $0<\de<\pi,\,\de\ne\pi/2,$ be an element in $A$.
We want to describe the set of billiard tables having caustics
$\Ga_{\de}$. The conjectures below aim at answering
Question~~\ref{number_quest1}.

\begin{conj}   \label{disjoint_roots_conj}
Let $m,n>1$ be distinct integers; let $S_m,S_n$ be the
corresponding polynomials in equation~~\eqref{filaseta_eq}. Then
their sets of complex roots are disjoint.
\end{conj}

The material of section~~\ref{trig_quant} and
section~~\ref{frac_lin_sub} yields  that
Conjecture~~\ref{disjoint_roots_conj} is equivalent to the
following claim.

\begin{conj}      \label{dist_tan_conj}
Let $m,n>1$ be distinct integers. Then equations $\tan mx=m\tan
x,\,\tan nx=n\tan x$ have no common solutions in $(0,\pi/2)$.
\end{conj}

For reader's convenience, we outline a proof that the two
conjectures are equivalent. Recall that $B_k$ denotes the set of
roots of the equation $\tan kx=k\tan x$ in $(0,\pi/2)$. By
Lemma~~\ref{prelim_lem}, Proposition~~\ref{polynom_prop}, and
Corollary~~\ref{trig_polynom_cor}, the set $\{\tan x:\,x\in B_k\}$
is the set of positive roots of the polynomial $R_k$. See
equation~~\eqref{R_n_eq}. Corollary~~\ref{gutkin_cor} provides a
fractional linear transformation that sends the positive roots of
$R_k,k>1,$ to the roots of $S_k$ in the semi-circle $\{|z|=1,\Im
z>0\}$. Now the information about the roots of $S_k$ contained in
Corollary~~\ref{filaseta_cor} implies the claim.

The following proposition lends support to
Conjecture~~\ref{dist_tan_conj}.

\begin{prop}      \label{small_k_prop}
Let $n>1$. We will say that a solution $x$ is nontrivial if $\tan
x\ne 0$.

\noindent {\em 1}. The system
\begin{equation}   \label{tan_n_n+k_eq}
\tan nx=n\tan x,\ \tan(n+k)x=(n+k)\tan x
\end{equation}
has no nontrivial solutions for $k=1,2$.

\noindent {\em 2}. The system
\begin{equation}   \label{tan_n_tan_k_eq}
\tan nx=n\tan x,\ \tan knx=kn\tan x
\end{equation}
has no nontrivial solutions for $k=2,3$.

\noindent {\em 3}. The systems
\begin{equation}   \label{tan_n_tan_2n+-1_eq}
\tan nx=n\tan x,\ \tan(2n\pm1)x=(2n\pm1)\tan x
\end{equation}
have no nontrivial solutions.
\begin{proof}
We have
$$
\tan(n+k)x=\frac{\tan nx + \tan kx}{1-\tan nx\tan kx}.
$$
Substituting this into equation~~\eqref{tan_n_n+k_eq} and using
that $\tan x\ne 0$, we obtain $1+n(n+1)\tan^2x=0$ in the case
$k=1$, and $(n+1)^2\tan^2x=0$ in the case $k=2$. This proves claim
1.

We have
$$
\tan 2nx=\frac{2\tan nx}{1-\tan^2nx},\ \tan 3nx=\frac{3\tan
nx-\tan^3nx}{1-3\tan^2nx}.
$$
Substituting these identities into
equation~~\eqref{tan_n_tan_k_eq}, and assuming $\tan x\ne 0$, we
obtain $1-n^2\tan x=1,\,8n^2\tan x=0$ if $k=2,\,k=3$ respectively.
This proves claim 2.

Equation~~\eqref{tan_n_tan_2n+-1_eq} and the identity
$$
\tan(2n\pm 1)x=\frac{\tan 2nx \pm \tan x}{1\mp\tan 2nx\tan x}
$$
yield the relationship $n^3\pm2n^2+n=0$ which has no solutions
$n>1$.
\end{proof}
\end{prop}

\begin{rem}      \label{3n+-1_rem}
{\em

A refinement of the above approach yields that the systems
%
%
$$
\tan nx=n\tan x,\ \tan(3n\pm1)x=(3n\pm1)\tan x
$$
%
do not have nontrivial solutions as well. The proof is rather
long, and we do not reproduce it here.

}
\end{rem}

Proposition~~\ref{small_k_prop} and Remark~~\ref{3n+-1_rem} yield
particular families of pairs of integers $m\ne n$ such that the
system $\tan mx=m\tan x,\,\tan nx=n\tan x$ has no nontrivial
solutions. This provides direct evidence supporting
Conjecture~~\ref{dist_tan_conj}. The work \cite{BFLT} provides
additional support for the equivalent
Conjecture~~\ref{disjoint_roots_conj}. We will now elaborate on
this.

Let $n\ge 4$. Set $\tS_n(x)=S_n(x)/(x-1)^3(x+1)$ if $n$ is odd and
$\tS_n(x)=S_n(x)/(x-1)^3$ if $n$ is even. By
Corollary~~\ref{filaseta_cor}, $\tS_n$ are polynomials with
integer coefficients; their roots are simple and belong to the
unit circle. Let $X$ be a property that holds for some natural
numbers. Denote by $\N(X)\subset\N$ be the set of natural numbers
having property $X$. We say that property $X$ {\em holds for
almost all positive integers} if $\N(X)\subset\N$ is a subset of
{\em density one}. A property that holds for {\em almost all pairs
of positive integers} is defined analogously.

The work \cite{BFLT} puts forward several conjectures about
irreducibility of polynomials over $\QQ$. It conjectures, in
particular, that polynomials $\tS_n$ are irreducible. See
Conjecture~~3 in \cite{BFLT}. Let $m\ne n$ be natural numbers. We
will say that {\em Conjecture~~\ref{disjoint_roots_conj} holds for
the pair $m,n$} if the sets of complex roots of the polynomials
$S_m,S_n$ are disjoint.

\begin{prop}       \label{alm_all_prop}
Conjecture~~\ref{disjoint_roots_conj} holds for almost all pairs
of positive integers.
\begin{proof}
By the preceding discussion, it suffices to show that for almost
all pairs $m\ne n$ the root sets of $\tS_m,\tS_n$ are disjoint.
Let $\II\subset\N$ be the set of integers $k$ such that $\tS_k$ is
irreducible. Let $\JJ\subset\II\times\II$ be the set of distinct
pairs. By Theorem~~4 in \cite{BFLT}, $\II\subset\N$ is a set of
density one. Thus, the sets
$\JJ\subset\II\times\II\subset\N\times\N$ have density one. But
for pairs $(m,n)\in\JJ$ the polynomials $\tS_m,\tS_n$ have
disjoint root sets.
\end{proof}
\end{prop}

Our next conjecture addresses Question~~\ref{number_quest2}.

\begin{conj}      \label{irrat_conj}
Let $x\ne 0$ satisfy $\tan nx = n\tan x$ for some $n>1$. Then $x$
is $x/\pi$ is irrational.
\end{conj}

The motivation for Conjecture~~\ref{irrat_conj} is as follows. If
$0<x<\pi/2$ satisfies equation~~\eqref{tan_nx_eq} then there is a
continuous family\footnote{A one-parameter family, if
Conjecture~~\ref{dist_tan_conj} holds.} of billiard tables $\Om$
with the invariant curve $\Ga_x$. The billiard map on $\Om$,
restricted to $\Ga_x$, is the rotation by $2x$. See the proof of
Proposition~~\ref{const_caust_prop}. If $x/\pi$ is rational, then
we acquire a lot of examples of billiard tables with invariant
circles realizing {\em rational rotations}. The billiard
literature indicates that such invariant circles are extremely
rare \cite{In88,GKn96}.


Let us now look at examples. The smallest $n$ for which equation
$\tan nx = n\tan x$ has nontrivial solutions is $4$. Below we
analyze its solutions for $n=4,5$.

\begin{exa}   \label{small_n_exa}
{\em

Set $z=\tan x$. From the recurrence relations
equation~~\eqref{recurr_eq}, we easily obtain
$P_4=4z-z^3,Q_4=1-6z^2+z^4,P_5=5z-10z^3+z^5,Q_5=1-10z+5z^4$. Thus,
the equation $\tan 4x=4\tan x$ is equivalent to $z^4-5z^2=0$.
Since $z\ne0$, this yields $z=\pm\sqrt{5}$. Denote by $x_4$ the
solution of $\tan 4x=4\tan x$ in $(0,\pi/2)$. Then
$x_4=\arctan(\sqrt{5})=\arcsin(\sqrt{5}/\sqrt{6})$. We analyze the
equation $\tan 5x=5\tan x$ is d the same way. Denote by $x_5$ the
unique solution of $\tan 5x=5\tan x$ in $(0,\pi/2)$. Then
$x_5=\arctan(\sqrt{5/3})=\arcsin(\sqrt{5}/2\sqrt{2})$.

Using the formula $\tan\frac{\pi}{5}=\sqrt{5-2\sqrt{5}}$, we
obtain the bounds
$$
\frac{\pi}{4}<x_5<\frac{3\pi}{10}<\frac{\pi}{3}<x_4<\frac{\pi}{2}.
$$
These inequalities do not imply that the numbers $x_4/\pi,x_5/\pi$
are  irrational; however, they show that the denominators cannot
be small.

}
\end{exa}

Our next proposition provides more evidence for
Conjecture~~\ref{irrat_conj}.

\begin{prop}    \label{pi_over_four_prop}
Let $\Om$ be a regular billiard table. Suppose that $\Ga_{\pi/4}$
or $\Ga_{\pi/3}$ is a caustic for $\Om$. Then $\Om$ is circular.
\begin{proof}
We will show that $\pi/3,\pi/4$ do not satisfy
equation~~\eqref{sin_eq} for any $n>1$. Set $\de=\pi/4$ and
examine both sides of equation~~\eqref{sin_eq} for $n=2,3,\dots$.
By periodicity of the sine function, everything is determined by
the residue $n\mod 4$. Let $n\equiv1\mod4$. Then in the left hand
side of equation~~\eqref{sin_eq} the numerator is $\sin(k\pi)=0$;
in the right hand side of equation~~\eqref{sin_eq} the numerator
is $\sin(k\pi+\pi/2)=\pm1$. Analogous considerations show that for
$n\equiv3\mod4$ equation~~\eqref{sin_eq} is not satisfied.

Let now $n\equiv0\mod4$. Thus $n-1\equiv3\mod4,n+1\equiv1\mod4$.
Then one of the two numerators in equation~~\eqref{sin_eq} is
$\pm1$ while the other is $\mp1$, hence equation~~\eqref{sin_eq}
does not hold. An analogous argument disposes of the case
$n\equiv2\mod4$. This proves the claim for $\pi/4$. The argument
for $\de=\pi/3$ follows the same pattern, and we leave it to the
reader.
\end{proof}
\end{prop}

\section{Conditional implications for the billiard and the floating}       \label{cond_impli}
We will now deduce some implications of the above conjectures to
billiard dynamics and capillary floating. We begin with the
billiard.

\subsection{Billiard tables with constant angle caustics}
\label{impli_billi_sub}
\hfill \break

\begin{prop}       \label{cond_const_caust_prop}
Let $\de\in A\setminus\{\pi/2\}$, and let $\Om\subset\RR$ be a
noncircular, regular billiard table with the caustic $\Ga_{\de}$.

\medskip

Suppose that Conjecture~~\ref{dist_tan_conj} holds. Then there is
a unique integer $n\ge 4$ and a unique parameter $0<\tau<1$ such
that $\Om$ is conformally equivalent to the table $\Om_{n,\tau}$
given by equations~~\eqref{table_eq1},~~\eqref{table_eq2}.
\begin{proof}
Let $\rho(\al)$ be the radius of curvature function for $\bo\Om$.
By Theorem~~\ref{const_caust_thm}, there are unique constants
$a,b,c$ satisfying $0<a^2+b^2<c$ such that $\rho(\al)=c+a\cos n\al
+ b\sin n\al$. The claim now follows from
Corollary~~\ref{1_param_prop}.
\end{proof}
\end{prop}

We say that a set $\Om\subset\RR$ is {\em rotationally symmetric
of order $n>1$} if  there is $P_0\in\RR$ such that $\Om$ is
invariant under the group $\Z/n\Z$ of rotations of $\RR$ about
$P_0$.

\begin{thm}     \label{cond_const_caust_thm}
Let $\Om\subset\RR$ be a noncircular, regular billiard table. The
following statements are equivalent.

\medskip

\noindent {\em 1}. The table $\Om$ has a caustic $\Ga_{\de}$,
$\de\ne\pi/2$.

\noindent {\em 2}. There is $n>3$ such that the Fourier
coefficients of the radius of curvature $\rho(\cdot)$ of $\bo\Om$
satisfy i) $c_n\ne 0$; ii) $c_k=0$ for all positive $k\ne n$.

\noindent {\em 3}. There is $n>3$ and $0<\tau<1$ such that $\Om$
is conformally equivalent to the table $\Om_{n,\tau}$ given by
equations~~\eqref{table_eq1},~~\eqref{table_eq2}.
\begin{proof}
Proposition~~\ref{cond_const_caust_prop} proves the implication
$1\Rightarrow 3$, while $2\Rightarrow 1$ is a byproduct of
Corollary~~\ref{1_param_prop}. The implication $3\Rightarrow 2$ is
obvious.
\end{proof}
\end{thm}

The preceding propositions rely on
Conjecture~~\ref{dist_tan_conj}.\footnote{Or the equivalent
Conjecture~~\ref{disjoint_roots_conj}.} The proof of the following
claim relies on Conjecture~~\ref{irrat_conj}.

\begin{thm}     \label{irrat_rot_thm}
Let $\Om\subset\RR$ be a regular billiard table, where $\bo\Om$ is
not a curve of constant width. Suppose that the table $\Om$ has a
caustic $\Ga_{\de}$ of constant type. Then the restriction of the
billiard map on $\Om$ to $\Ga_{\de}$ is an irrational rotation.
\begin{proof}
By the proof of Proposition~~\ref{const_caust_prop}, the
transformation in question is the rotation by $2\de/\pi$. The
claim now follows from Conjecture~~\ref{irrat_conj}.
\end{proof}
\end{thm}

From now until the end of the section, we will assume the truth of
all conjectures in section~~\ref{conj_sub}.

\begin{thm}     \label{conseq_thm}
There is a dense countable set $R\subset(0,1)$ of irrational
numbers such that the following claims hold.

\medskip

\noindent {\em 1}. For every $\rho\in R$ there is a one-parameter
family of regular billiard tables $\Om_{\tau}$ having a constant
angle caustic with the rotation number $\rho$. The curves
$\bo\Om_{\tau}$ are real analytic. Every regular billiard table
having a constant angle caustic with the rotation number $\rho$ is
conformally equivalent to a unique table $\Om_{\tau}$.

\noindent {\em 2}. Let $\rho\in(0,1)\setminus R$. Suppose that a
regular billiard table $\Om$ has a constant angle caustic with the
rotation number $\rho$. i) If $\rho=1/2$ then $\Om$ has constant
width. ii) If $\rho\ne 1/2$ then $\Om$ is a disc.
\begin{proof}
Claim 1 is immediate from Theorem~~\ref{cond_const_caust_thm} and
Theorem~~\ref{irrat_rot_thm}. Claim 2 follows by combining these
statements with Theorem~~\ref{const_caust_thm}.
\end{proof}
\end{thm}
\subsection{Two-dimensional capillary floating}
\label{impli_float_sub}
\hfill \break In what follows we assume that all of the
conjectures in section~~\ref{conj_sub} hold.

\begin{thm}    \label{float_thm}
Let $\Om\subset\RR$ be a regular, compact, convex domain. Then the
following holds.

\medskip

\noindent 1. Suppose that $\Om$ is not a disc. Then $\Om$ floats
in neutral equilibrium at any orientation with the contact angle
$\pi-\de\ne\pi/2$ if and only if  there is $n>3$ and $0<\tau<1$
such that $\Om$ is conformally equivalent to the domain
$\Om_{n,\tau}$ given by
equations~~\eqref{table_eq1},~~\eqref{table_eq2}.

\medskip

\noindent 2. Suppose that $\Om$ is not a domain of constant width.
If $\Om$ floats in neutral equilibrium at any orientation  with
the contact angle $\ga$ then $\ga/\pi$ is irrational.

\medskip

\noindent 3. There is a countable dense set $A\subset(0,\pi)$
containing $\pi/2$ and symmetric about this point such that the
following holds:

\medskip

\noindent If $\Om$ floats in neutral equilibrium at any
orientation with the contact angle $\ga\in(0,\pi)\setminus A$ then
$\Om$ is a disc.
\begin{proof}
The claims are the counterparts of statements in
Theorems~~\ref{cond_const_caust_thm},~~\ref{irrat_rot_thm},
and~~\ref{conseq_thm}.
\end{proof}
\end{thm}

\vspace{5mm}

\medskip

\end{document}